\documentclass[onefignum,onetabnum,a4paper,twocolumn]{article}
\usepackage{epsfig}
\usepackage{epstopdf}
\usepackage{amsmath,amssymb}
\usepackage{textcomp}
\usepackage{gensymb}

\title{Stabilising millennial oscillations in large-scale ocean circulation with a delayed feedback due to a circumpolar current} 

\author{{\sc Andrew Keane$^{1,2}$\footnote{corresponding author:  {\tt andrew.keane@ucc.ie}}, Alexandre Pohl$^{3}$, Henk A. Dijkstra$^{4,5}$ and Andy Ridgwell$^{6}$} \\ 
$^{1}$School of Mathematical Sciences, University College Cork, Cork, Ireland\\
$^{2}$Environmental Research Institute, University College Cork, Cork, Ireland\\
$^{3}$Biog\'eosciences UMR 6282, Universit\'e Bourgogne Europe, CNRS, F-21000, Dijon, France\\
$^{4}$Institute for Marine and Atmospheric research Utrecht, Department of Physics\\Utrecht University, Utrecht, The Netherlands\\
$^{5}$Center for Complex Systems Studies, Department of Physics, Utrecht University\\Utrecht, The Netherlands\\
$^{6}$Department of Earth and Planetary Sciences, University of California, Riverside, CA, USA\\
}

\date{April 2025}

\begin{document}

\maketitle

\begin{abstract}
The global ocean circulation plays a pivotal role in the regulation of the Earth's climate.
The specific pattern and strength of circulation also determines how carbon and nutrients are cycled and via the resulting distribution of dissolved oxygen, where habitats suitable for marine animals occur.
However, evidence from both geological data and models suggests that state transitions in circulation patterns  have occurred in the past. Understanding the controls on marine environmental conditions and biodiversity requires a full appreciation of the nature and drivers of such  transitions.
Here we present stable millennial oscillations of meridional overturning circulation in an Earth system model of intermediate complexity, cGENIE, that appear to only occur in the presence of a circumpolar current.
To demonstrate that a circumpolar current can act as a driver of stable oscillations, we adapt a simple ocean box model to include a delayed feedback to represent the effect of a circumpolar current on meridional overturning circulation. We investigate the millennial oscillatory solutions that arise in the box model by bifurcation analysis and show that the model can reproduce the same bifurcation structure observed in the Earth system model.
Our results provide new insights into the nature of oscillations that could have occurred under certain continental configurations  in the geological past, and also highlight the potential influence of the changing Antarctic circumpolar current speed on the stability of the Atlantic meridional overturning circulation.
\end{abstract}




\section{Introduction}
\label{section:intro}

The ocean circulation plays a critical role in the Earth's global climate system. Wind-driven currents acting over the first hundreds of metres impose the first-order pattern and magnitude of ocean heat transport, and therefore largely determine superficial climate and notably the latitudinal gradient of temperature \cite{Rose2013}. Upper-ocean currents also control equatorial and coastal upwelling position and strength, thus the spatial patterns of primary productivity and extension of oxygen minimum zones \cite{Pohl2017}. 
The deep ocean circulation, on the other hand, is affected by density differences and is important for the redistribution of nutrients in the global ocean (and critically, their return back to the ocean surface) \cite{pohl22} as well as the sequestration of atmospheric carbon in the ocean interior \cite{Bouttes2010}.
Overturning circulation also controls the transport of oxygen equilibrated in surface seawater into the ocean interior (`ventilation') and, ultimately, the capacity of the ocean  to host animal life  at depth \cite{Pohl2021,Stockey2021}. 

Abundant sedimentological, geochemical, and paleontological data demonstrate that the global ocean circulation has significantly varied through geological time. At the (multi-)million-year scale, it has varied in response to continental rearrangement. For instance, Trabucho Alexandre \emph{et al.} \cite{TrabuchoAlexandre2010} propose that intense deoxygenation and the deposition of organic-rich sediments developed in the mid-Cretaceous North Atlantic, around 100 million years ago, as a consequence of a specific ocean circulation pattern making this ocean basin act as a nutrient trap. A similar mechanism was suggested by Meyer \emph{et al.} \cite{Meyer2008a} to explain photic-zone euxinia expansion in the Paleo-Tethys ocean at the Permian-Triassic boundary. At a higher frequency, the ocean circulation, both shallow and deep, is also known to vary in concert with global climate. Today, vigorous deep-water formation takes place in the North Atlantic (and Southern) Ocean while limited ventilation is observed at similar northern latitudes in the  Pacific Ocean \cite{Rae2020}. During the Last Glacial Maximum, sediment core proxy data for ocean ventilation suggest that the North Pacific Ocean was better ventilated at intermediate depths \cite{Rae2020}, while the Atlantic circulation was probably more sluggish (shallower and/or less vigorous) \cite{Howe2016}. Similarly, Freeman \emph{et al.} \cite{Rae2020} used radiocarbon measurements to reconstruct the benthic--planktonic $^{14}$C age offset (i.e. `ventilation age') 
during the last deglaciation and demonstrated a see-saw in the ventilation of the intermediate Atlantic and Pacific Oceans. 

In order to assimilate the proxy database in a quantitative framework and gain mechanistic understanding of the drivers of observed changes in ocean properties, (paleo)climatologists use ocean-atmosphere general circulation models (GCMs) \cite{Donnadieu2016,Ferreira2011,Laugie2021}. These models provide complex representations of ocean physics and dynamics and the coupling with the atmospheric component, sea ice, and possibly ocean biogeochemistry and land surface \cite{Sepulchre20209}. Models of the Coupled Model Intercomparison Project phase 6 (CMIP6), for instance, aim to offer a satisfactory representation of the climate system \cite{Zelinka2020}. However, these models are computationally expensive and consequently, are generally run only for a few hundred to a thousand years. Because this duration approaches the ocean mixing timescale, it is assumed that it allows the model to approximately reach deep-ocean equilibrium, i.e., the physical characteristics of the whole ocean (temperature and salinity), and thus ocean dynamics, are supposed stable from that point in time \cite{Lunt2016,Donnadieu2016}. By doing so, it is also implicitly assumed that the ocean circulation does in fact eventually reach a stable steady state. Yet, various attempts to run global climate models for a longer time period demonstrated more complex behaviours of the ocean circulation, which sometimes exhibit regimes of stable self-sustained oscillations. 
Oscillations have notably been simulated in ocean-only GCMs \cite{Sirkes2001,Weijer2003} and ocean GCMs coupled to simple atmospheric components (Earth system models of intermediate complexity; EMICs) \cite{Meissner2008,Haarsma2001}. 
To our knowledge, the only studies reporting such self-sustained oscillations based on ocean-atmosphere GCM simulations are Refs.~\cite{armstrong22,valdes21} using the fully coupled HadCM3 GCM, Ref.~\cite{kuniyoshi2022effect} using the MIROC4m GCM, Ref.~\cite{klockmann2020coupling} using the Max Planck Institute Earth system model and the last glacial CESM simulations of Vettoretti \emph{et al.} \cite{Vettoretti2022}. 
The absence of further documentation of stable oscillatory regimes in other fully coupled GCMs is not surprising considering that a simulation duration of at least 6 thousand years seems to be required to identify stable oscillations, while a typical run-time for GCMs used to explore past (geological) climate states is typically only about 3 thousand years (e.g. \cite{chen22,Laugie2021,wong21}). 

The existence of regimes of stable self-sustained oscillations in ocean circulation has recently been demonstrated by Pohl \emph{et al.} \cite{pohl22} for specific combinations of past continental configuration and climatic state using the EMIC cGENIE. The oscillations have been invoked to explain both the poorly-oxygenated oceans that pervaded several hundreds of millions of years ago \cite{Tostevin2020}, when marine life started diversifying in the marine realm \cite{Alroy2008}, as well as enigmatic thousand-year scale cycles identified in marine sediments of the deep geological past \cite{Dahl2019}. Based on recent studies that pinpoint the tight coupling between marine redox conditions and marine biodiversity \cite{Penn2018,Stockey2021,Deutsch2015a,Rubalcaba2020}, the authors conclude that the simulated high-frequency variability in ocean redox conditions would have constituted harsh conditions for the development of marine life at that time. The study by Pohl \emph{et al.} \cite{pohl22} exemplifies the implications that such oscillations have for our understanding of the evolution of both global climate and marine biodiversity, but also highlights the need to better characterise the physical mechanisms causing these oscillations. It focused on the consequences of the resulting state transitions in ocean circulation for ocean oxygenation and the interpretation of geological redox indicators, rather than on the underlying mechanisms of the oscillations. In addition, the experimental setup used by the authors, involving complex land mass shapes and configurations, as well as the presence in the model of a variety of  climatic and biogeochemical feedbacks, creates challenges for a detailed analysis. As a result, the physical mechanisms causing the stable self-sustained oscillations and state transitions in ocean dynamics in cGENIE remain to be fully elucidated.

Here, we analyse the mechanisms at play by combining a series of new Earth system model experiments conducted using an idealised ocean configuration, with a box model designed specifically to reproduce and identify key characteristics of the cGENIE simulations. Experiments using cGENIE with the idealised continental configuration  serve as a stepping-stone between results of oscillations in simulations using ``real-world'' continental configurations \cite{pohl22} and understanding gained through the analysis of simple ocean box models. 

Box models have been often used to identify physical mechanisms 
of oscillatory behaviour in the ocean-climate system. Regarding variability of the  meridional overturning circulation, two basic types of oscillatory behaviour  have  been identified. Overturning (or loop) oscillations can be found  in a four-box model \cite{Tziperman1994} and  involve the advective propagation of a buoyancy anomaly along the overturning loop. Such variability is of centennial time scale  when only the Atlantic is considered and millennial scale when the global overturning is considered \cite{Weijer2003}. Oscillatory variability  may also be induced by convective processes. The basics of the behaviour can be understood from the so-called `flip-flop' oscillation \cite{Welander1982}, where oscillations between convective 
and non-convective states occur. The origin of the millennial time-scale oscillations,  also referred to as `deep-decoupling' oscillations or `flushes' is, however, in more detailed models best  explained by  a box model which combines both advective and convective processes \cite{ColindeVerdiere2007}. 

Here, we present a model that reproduces the oscillations obtained in the cGENIE Earth system model \cite{pohl22} in a very simple way, by embedding key climatic features in a simple 3-box model. The novel element is that a time delay is introduced, which serves as an efficient representation of the effects of a circumpolar current on the oceanic salt transport.

\section{EMIC experiments}
\label{section:genie}

\begin{figure}[t!]
  \centering
  \vspace*{1mm}
  \includegraphics[width=\columnwidth]{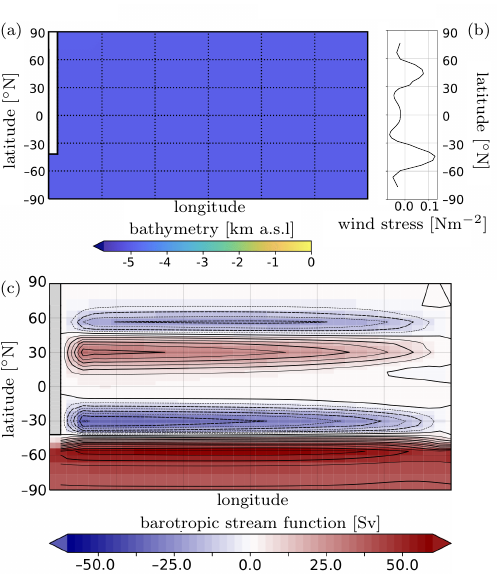}
  \caption{cGENIE Drake World setup. (a) Model bathymetry. Landmasses (above sea level; a.s.l.) are shaded white. (b) Zonal wind stress field used as boundary condition. 
  (c) Barotropic streamfunction (Sverdrup, 1 Sv = 10$^6$ m$^{-3}$ s$^{-1}$) simulated at 16 times the pre-industrial atmospheric CO$_2$ concentration of 280 ppm. The boundary is periodic in the horizontal direction.
}
  \label{fig:DW}
\end{figure}

To gain insight into the existence of self-sustained oscillations in the ocean circulation, numerical simulations are conducted using an idealised setup of the EMIC cGENIE \cite{Ridgwell2007}. cGENIE consists of a reduced physics (frictional geostrophic) 3D ocean circulation model coupled to a 2D energy-moisture-balance atmospheric component plus sea-ice module. Although we do not show the results of it here, cGENIE also includes a state-of-the-art representation of marine biogeochemistry. Compared to CMIP6 Earth system models, the simplified climatic component of cGENIE offers a rapid model integration time, permitting the model to be run for tens of thousands of years.

The adopted idealised cGENIE setup is ``Drake World''. This differs slightly from the classic MIT GCM implementation \cite{Ferreira2010}, but is primarily characterised by a thin strip of land extending southwards to 45\degree S and a 5-km deep flat-bottomed ocean everywhere else (Fig.~\ref{fig:DW}(a)). The model uses a 36$\times$36 equal-area horizontal grid, with 16 unevenly-spaced $z$-coordinate vertical levels in the ocean. We apply idealised boundary conditions of zonally averaged wind stress and speed, plus a zonally averaged planetary albedo, following Vervoort \emph{et al.} \cite{Veervort2024} and van de Velde \emph{et al.} \cite{Velde2021} (Fig.~\ref{fig:DW}(b)). The solar constant is set to the modern value (1368 W m$^{-2}$). The physical parameters controlling the model climatology  follow Cao \emph{et al.} \cite{Cao2009}. Additional experiments are conducted using an alternative, idealised continental configuration with a strip of land extending from the North Pole to the South Pole (``Ridge World'' \cite{Ferreira2010, pohl22}) and everything else unchanged.

Figure~\ref{fig:DW}(c) shows the barotropic stream-function (units of Sverdrup, 1 Sv = 10$^6$ m$^{3}$ s$^{-1}$), as a measure of the vertically averaged flow, in Drake World with atmospheric CO$_2$ at 16 times the pre-industrial atmospheric level (PAL, 1 PAL $= 280$ ppm).
As with all the following results, this shows the state of the model following spin-up (i.e. after all transient behaviour has died away).
The most notable feature of the barotropic stream-function is the strong zonal flow in the southern channel, similar to the Antarctic Circumpolar Current observed in the Earth's modern continental configuration. It is a mainly wind-driven current that constantly transports water around the pole and persists also for the other levels of atmospheric CO$_2$ considered here.

\begin{figure}[t!]
  \centering
  \includegraphics[width=\columnwidth]{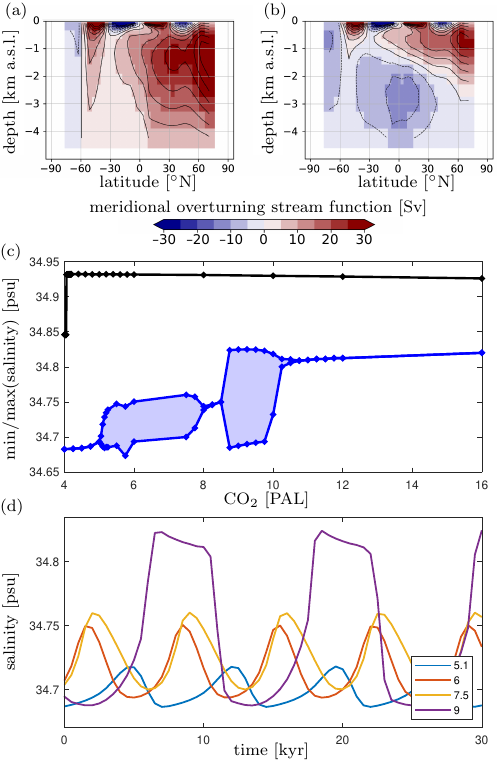}
  \caption{cGENIE results. 
  Meridional overturning streamfunction for Drake World at (a) 4 PAL and (b) 16 PAL. (c) Envelope of mean salinity (across the southern channel: comprising all depths, longitudes and 45--90\degree S) simulated at various atmospheric CO$_2$ levels, illustrating regimes of stable equilibria (where no envelope is visible) and stable oscillations (where the envelope can be seen), for Drake World (blue) and Ridge World (black). (d) Example time series of stable oscillations (legend indicates PAL) simulated using Drake World with transients (initial spin-up of 30 kyr) removed.
}
  \label{fig:genie}
\end{figure}

We find that when atmospheric CO$_2$ forcing is varied in the model (from 4 to 16 PAL), ocean dynamics exhibits a complex response. 
Overall, we observe a gradual shift in the locus of deep-water formation from the (oceanic) North Pole (Fig.~\ref{fig:genie}(a)) to the South Pole (Fig.~\ref{fig:genie}(b)). 
Interestingly, during this shift we encounter two families of stable oscillations.

The blue points in Fig.~\ref{fig:genie}(c) show mean salinity for an ensemble of simulations in Drake World. Throughout this paper, we take the mean salinity across the southern channel, comprising all depths, longitudes and 45--90\degree S.
Because the simulations take a long time to run (around 6 days for a 60-kyr long cGENIE simulation as shown in Fig.~\ref{fig:genie}(c)), they are run in parallel for different atmospheric CO$_2$ levels using the same initial conditions.
As CO$_2$ increases, the ocean transitions from a stable steady state with northern deep-water formation to a regime of stable oscillations between northern deep-water formation and a weakened circulation state (represented by the envelopes in Fig.~\ref{fig:genie}(c)). This weakened state becomes stable near 8.5 PAL. Further increasing of CO$_2$ gives rise to another regime of stable oscillations between a largely collapsed circulation state and strong circulation with southern deep-water formation. Finally, the southern deep-water formation steady state becomes stable for higher levels of atmospheric CO$_2$.
Example time series of the oscillations are shown in Fig.~\ref{fig:genie}(d).
Although results are only shown here for salinity across the southern channel, these changes in ocean circulation have been shown to induce significant and global periodic changes in many other observables, such as deep-ocean ventilation and oxygenation (see Ref.~\cite{pohl22}, extended data Fig.~10).

An additional series of experiments was conducted on Ridge World, i.e., under identical boundary conditions but with the strip of land extending to the southernmost grid points. The results reveal a transition from a stable steady state with northern deep-water formation to one with southern deep-water formation. This transition is sudden, represented by the jump of black points in Fig.~\ref{fig:genie}(c) near 4 PAL, and is akin to the classic ``tipping point'' behaviour of the Atlantic meridional overturning circulation observed in many studies \cite{rahmstorf1995bifurcations, alkhayuon2019basin, van2024physics}. However, importantly, this transition does not exhibit regimes of stable oscillations, at least not for the levels of CO$_2$ tested. A primary question addressed in this paper is then: Why do these oscillations appear once the ocean has an open-ocean passage around the South Pole?

\section{Simple box model}
\label{section:model}

We begin with the simplified ocean box model from section~3 of Ref.~\cite{titz02}, which consists of two polar boxes and one upper equatorial box, as illustrated below in Fig.~\ref{fig:boxes}. Assuming salt conservation, the model is written in terms of the salt content of the southern polar box (box 1) and the northern polar box (box 2).
The model, as it appears in Ref.~\cite{titz02}, represents only meridional transport and its behaviour corresponds to the meridional overturning circulation of Ridge World. In Drake World, on the other hand, the presence of a circumpolar current leads to the zonal transport of salt in the southern channel.

Here, we adapt the ocean box model to capture the effect of this zonal transport by including a delayed feedback term on the salinity of box~$1$. In particular, a parcel of water in the southern polar box will leave its position at time $t$, travel around the southern pole, and return to the same position at time $t+\tau$. 
This creates a delayed effect on the average salinity in box~$1$.
The model becomes a set of delay differential equations (DDEs) for the salinity variables $S_{i\in \{1,2\}}$:
\begin{align}
\label{eq:titz_delay}
V\dot{S_1} &= S_0 F_1 + m \left(S_2 - S_1\right) + \sigma \left(S_1^{\tau} - S_1\right)\\ \nonumber
V\dot{S_2} & =  -S_0 F_2 + m \left(S_3 - S_2\right).
\end{align}
Parameter $V$ is the (equal) volume of the boxes and $S_0$ is a reference salinity. 
Parameters $F_1$ and $F_2$ represent surface freshwater fluxes between the boxes due to atmospheric water vapour transport and wind-driven surface currents. The delayed variable $S_1^{\tau} = S_1\left(t-\tau\right)$ has delay time $\tau$ and $\sigma$ is the feedback strength. 
In other words, we have water with salt density $S_1$ leaving and salt density $S_1^{\tau}$ entering box~1 with $\sigma$ representing zonal flux.
The variable $m$ represents the volume transport between boxes due to the density gradient across the ocean:
\begin{equation}
\label{eq:m}
    m = k\left[\beta\left(S_2 - S_1\right) - \alpha T^*\right].
\end{equation}
Parameter $k$ is a hydraulic constant, while $\beta$ and $\alpha$ are salinity and temperature expansion coefficients, respectively. The model assumes a fixed temperature gradient $T^*$ between the northern and southern boxes. 
In Fig.~\ref{fig:boxes}, a positive $m$ corresponds to clockwise circulation. Note that here we only discuss the positive $m$ solutions, and therefore consign the negative $m$ equations to \ref{app:neg_m} for completeness.
Assuming that salt is conserved, we have $S_3 = 3S_0 - S_1 - S_2$; see \ref{app:salt} for a discussion of the validity of this salt conservation rule. 

\begin{figure}[t]
  \centering
  \vspace*{1mm}
  \includegraphics[width=\columnwidth]{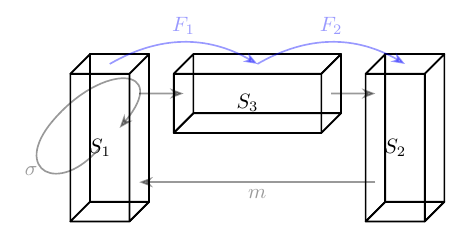}
  \caption{Boxes representing regions of the ocean with salinity variables $S_{i\in\{1,2,3\}}$ with arrows representing meridional water flux of strength $m$ and zonal flux of the southern channel of strength $\sigma$. Parameters $F_1$ and $F_2$ represent freshwater fluxes.
}
  \label{fig:boxes}
\end{figure}

In this paper we take parameter values from Ref.~\cite{titz02} and let $k=23\times 10^{17} \text{~m}^3\text{~yr}^{-1}$, $\alpha=1.7\times 10^{-4} \text{~K}^{-1}$, $\beta=0.8\times 10^{-3} \text{~psu}^{-1}$, and $S_0=35 \text{~psu}$. We let $V=3.5\times 10^{17} \text{~m}^3$, which is approximately the volume of the southern channel in Drake World. Furthermore, we set $F_2=1 \text{~Sv}$ ($10^6 \text{~m}^3\text{~s}^{-1}$) and $T^*=0 \text{~K}$, and note that qualitatively similar results to those presented below are found when varying these parameters. Finally, we approximate delay times for each latitude within the southern channel using averaged zonal ocean speeds from cGENIE simulations across longitude and depth. The average of these delay times gives us a rough estimate for $\tau$ of about $1000$ years. Therefore, we use similar values for $\tau$ in the results presented below. 

The presence of a circumpolar current creates a delayed effect on the average salinity in the southern polar region.
We acknowledge that a delayed feedback is not the only way to model the circumpolar current. For example, one could instead use multiple boxes to represent zonal flow across the southern channel. We suggest that a delayed feedback is a relatively elegant modelling approach, since it is mathematically concise, requiring only one additional term to the model, and allows for some simple physical interpretation. In the case of zero wind speed across the southern channel, there is no significant zonal flux, such that $\sigma$ approaches zero and the delayed feedback effect becomes negligible. 
In the case of $\sigma$ equal to zero, the model reverts to representing only meridional flow and corresponds to the behaviour of Ridge World rather than Drake World.
For non-zero wind speed across the southern channel, the wind speed drives the zonal flow and is, therefore, inversely proportional to the delay time $\tau$. 
The delayed feedback term does not change the steady-state solutions of the system, since $S_1^{\tau}=S_1$ for all steady states, such that the term vanishes. It may, nonetheless, affect the stability of a steady state. When $\tau$ is sufficiently small the delayed feedback may have a stabilising effect. For example, if $S_1(t)$ experiences a positive perturbation away from its steady-state value, then $S_1>S^{\tau}_1$ such that the delayed feedback works to decrease $S_1$, since $\dot{S}_1\propto \left(S^{\tau}_1-S_1\right)$. However, if $\tau$ and $\sigma$ are sufficiently large, this restoring behaviour of the delayed feedback may overshoot the steady state and instead promote oscillatory dynamics. This type of phenomenon is typical of delay-induced dynamics; see, for example, chapter~2 of Ref.~\cite{smith2011introduction}.
Of course, we could use additional delayed feedback terms to represent the influence of further periodic flows, such as the surface gyres seen in Fig.~\ref{fig:DW}(c). However, these surface effects would translate to delayed feedback terms with very small $\sigma$ and $\tau$ values, resulting in negligible effects for the meridional overturning, and are, therefore, neglected. It is only the periodic flow in the southern channel that, due to the lack of land mass, possesses enough momentum to generate a deep flow with sufficiently long delay times, resulting in large $\tau$ values.

Clearly, the box model representation is very much highly simplified compared to a dynamic 3D ocean circulation model. As such, we do not attempt to find quantitative agreement between the two models. Rather, we demonstrate that this simple model possesses surprisingly rich dynamics and is sufficient for a qualitative reproduction of the interesting ocean circulation behaviour observed in cGENIE experiments.

In Ref.~\cite{titz02} the model without the delayed feedback was shown to exhibit, not only the typical bistability between upper ($m>0$) and lower ($m<0$) branches of solutions, but also a Bogdanov-Takens bifurcation. This means that instead of the upper branch losing stability at a fold bifurcation, it may for certain values of freshwater flux already lose stability at a Hopf bifurcation, which is always subcritical \cite{titz02b}. The resulting unstable periodic orbits terminate at homoclinic bifurcations. We now investigate the effect that the delayed feedback has on this bifurcation structure. 
As shown in \ref{app:steady}, the analysis of model~(\ref{eq:titz_delay}) is complicated by the delay term, such that we primarily rely on numerical methods.
In the following section, solutions to the model and their stability properties are found using the continuation package DDE-Biftool~\cite{engelborghs00,janssens10,sieber14,wage14}.

\section{Stable oscillations}
\label{section:oscillations}

\begin{figure}[t!]
  \centering
  \vspace*{1mm}
  \includegraphics[width=\columnwidth]{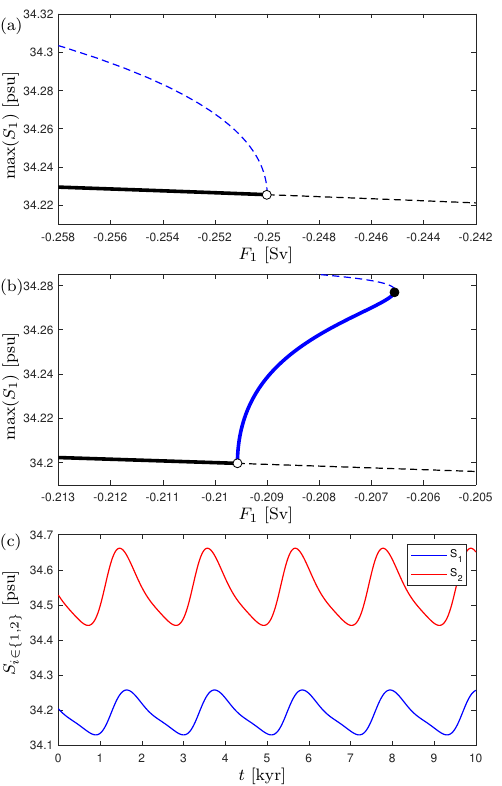}
  \caption{Solutions of model~(\ref{eq:titz_delay}) depending on freshwater forcing parameter $F_1$. Solid/dashed curves represent stable/unstable solutions. Black are equilibria and blue are maxima of periodic solutions. White and black filled circles represent Hopf and fold bifurcations, respectively. (a) is without the delayed feedback (i.e. $\sigma=0$), while (b) is with $\sigma=11 \text{~Sv}$. The delay time $\tau$ is set to $900$ years. (c) shows an example time series with $\sigma=11\text{~Sv}$ and $F_1=-0.208\text{~Sv}$.
}
  \label{fig:bif_F1}
\end{figure}

The solutions to the box model shown in Fig.~\ref{fig:bif_F1} demonstrate how the delayed-feedback effect allows the existence of stable periodic orbits. Black curves show the steady-state solutions of $S_1$ as parameter $F_1$ is varied, while the blue curves indicate the maximum of $S_1$ for periodic solutions. Solid and dashed curves indicate stable and unstable solutions, respectively. In panel~(a) the parameter $\sigma$ is set to zero; that is, there is no delayed-feedback effect. In this case, as $F_1$ increases, the stable steady-state solution loses stability at a subcritical Hopf bifurcation, resulting in a branch of unstable periodic orbits. When the delayed feedback is active, the Hopf bifurcation may become supercritical, resulting in a branch of stable periodic orbits, as shown in panel~(b). In this particular example, the stable periodic solutions lose stability at a fold of periodic orbits. The periods of the resulting stable oscillations have a millennial timescale as demonstrated by the example time series shown in panel~(c). 

Some readers familiar with control theory may not be surprised by the ability of the delayed feedback to stabilise the pre-existing periodic orbits. Curiously, the additional term in our model that represents the zonal flow takes the same mathematical form as Pyragas control. This control method is used to stabilise unstable periodic orbits, as has been demonstrated in many applications, for example, in semiconductor laser experiments \cite{pyragas92,schikora11}. 
To be precise, the delayed feedback effect in our model is not equivalent to Pyragas control. In the latter, the control is applied to all variables of the system (we only apply it to the salinity in one box) and the delay time is set equal to the period of the unstable periodic orbit to be stabilised. Nonetheless, the method may still be successful for other delay times, as investigated in Ref.~\cite{purewal14}. In model~(\ref{eq:titz_delay}) we find that the periods of the stable oscillations are generally at least twice the delay time.

\begin{figure}[t]
  \centering
  \vspace*{1mm}
  \includegraphics[width=\columnwidth]{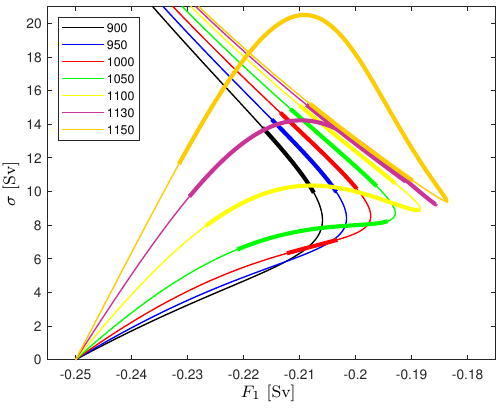}
  \caption{Curves of Hopf bifurcations. At a Hopf bifurcation the steady-state solution loses (or gains) stability and produces a branch of periodic solutions. Where the curves are thin/thick the Hopf bifurcation produces a branch of unstable/stable periodic solutions. The colours refer to different delay times $\tau$ shown in the legend. 
  }
  \label{fig:hopf}
\end{figure}

Whether the delayed feedback induces a change in criticality of the Hopf bifurcation depends on the feedback parameters $\sigma$ and $\tau$. Figure~\ref{fig:hopf} illustrates this dependence, where curves of Hopf bifurcations are calculated in the $(F_1,\sigma)$-plane for different values of $\tau$. The thick parts of the curves indicate that the Hopf bifurcation is supercritical. Notice that, depending on the feedback strength and delay time, it is possible that the system experiences multiple Hopf bifurcations as $F_1$ increases.

\begin{figure}[t!]
  \centering
  \vspace*{1mm}
  \includegraphics[width=\columnwidth]{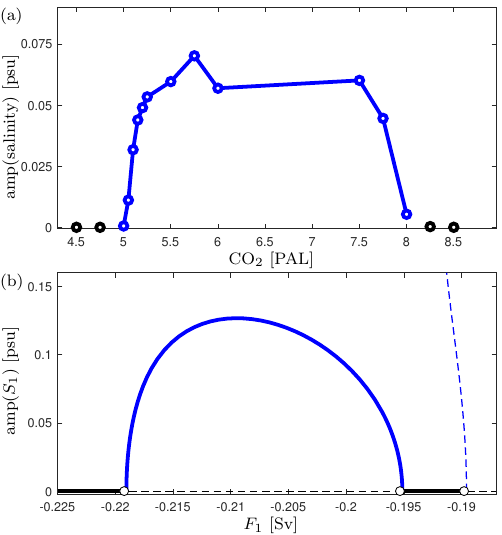}
  \caption{Amplitudes of oscillations with steady-state and periodic solutions shown in black and blue, respectively. (a) cGENIE experiments with varying atmospheric CO$_2$. (b) Model~(\ref{eq:titz_delay}) solutions while varying $F_1$ with $\sigma=9.5 \text{~Sv}$ and $\tau=1100\text{~years}$. Solid/dashed curves represent stable/unstable solutions. White filled circles represent Hopf bifurcations.
}
  \label{fig:amp_small}
\end{figure}

Figure~\ref{fig:amp_small}(a) is a bifurcation diagram, or more precisely an ``attractor diagram'', showing the amplitude of oscillations in terms of mean salinity (across the southern channel of Drake World) constructed from simulation results of cGENIE. 
We observe that as CO$_2$ is increased the steady-state solution loses stability near 5 PAL. The amplitude growth near this bifurcation point resembles a square-root law, typical of a Hopf bifurcation \cite{strogatz18}. The family of periodic solutions disappears at, what appears to be, another supercritical Hopf bifurcation near 8 PAL. The bifurcation point must be very close to 8 PAL, since it was necessary to run the 8 PAL simulation for an extremely long time (up to 120 thousand years) while the amplitude of the oscillations continued to decrease very slowly. Finally, the steady-state solution again loses stability after 8.5 PAL.

As illustrated in Fig.~\ref{fig:hopf} we may choose parameter values so that the model will pass through two supercritical Hopf bifurcations as $F_1$ increases. Realistically, $F_1$ is not the only parameter value that would change due to changes in atmospheric CO$_2$. Parameters $F_2$, $T^*$, $\sigma$ and $\tau$ would all be susceptible to influence from atmospheric changes. Nonetheless, we only consider varying $F_1$ in an effort to keep the analysis simple and demonstrate qualitative changes in behaviour. Figure~\ref{fig:amp_small}(b) shows that the effect of the circumpolar current is sufficient for reproducing a family of periodic solutions between two supercritical Hopf bifurcations. In the box model the steady-state solution finally loses stability at a subcritical Hopf bifurcation.

\begin{figure}[t!]
  \centering
  \vspace*{1mm}
  \includegraphics[width=\columnwidth]{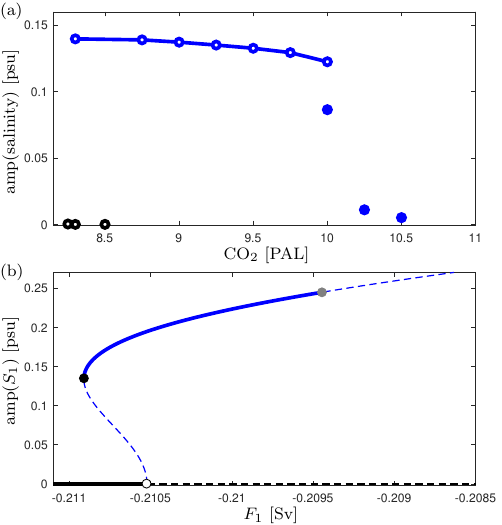}
  \caption{Amplitudes of oscillations with steady-state and periodic solutions shown in black and blue, respectively. (a) cGENIE experiments with varying atmospheric CO$_2$ with steady-state and periodic solutions shown as black and blue circles, respectively. Filled blue circles indicate solutions not captured by the box model. (b) Model~(\ref{eq:titz_delay}) solutions while varying $F_1$ with $\sigma=9 \text{~Sv}$ and $\tau=850\text{~years}$. Solid/dashed curves represent stable/unstable solutions. White, black and grey filled circles represent Hopf, fold and period-doubling bifurcations, respectively.
}
  \label{fig:amp_large}
\end{figure}

Figure~\ref{fig:amp_large}(a) is the same as Fig.~\ref{fig:amp_small}(a), only for higher atmospheric CO$_2$. Near 8.5 PAL we observe a jump to large-amplitude oscillations, as well as bistability between steady-state and oscillatory solutions. Both the jump to large-amplitude oscillations and bistability are reproduced by the box model in panel~(b). Here, the steady-state solution loses stability in a subcritical Hopf bifurcation. The unstable periodic orbits become stable via a fold bifurcation of periodic orbits, before later losing stability at a period-doubling bifurcation.

The termination of the large-amplitude oscillations in Fig.~\ref{fig:amp_large}(a) involves co-existing oscillatory solutions at 10 PAL and a gradual transition to a state of strong convection at the southern pole. These solutions are shown as filled blue circles in panel~(a). Following these oscillatory solutions, for larger values of CO$_2$, we only find steady-state solutions that correspond to $m<0$ solutions in the box model, which we do not consider here for simplicity. The re-creation of transitions between strong and weak convecting states in the box model would require upper and lower boxes, as well as the introduction of an additional convection timescale into the dynamics as in Ref.~\cite{ColindeVerdiere2007}, and is therefore beyond the scope of model~(\ref{eq:titz_delay}). 

\section{Discussion}
\label{section:disc}

The delayed feedback provides the simple box model with the capacity to mimic important qualitative behaviour observed in cGENIE experiments using an idealised configuration. It hence provides a simple mechanism for the occurrence of stable oscillations in the ocean circulation and its dependence on a zonal flow across the southern channel. The analysis, conducted using the numerical continuation software DDE-Biftool, highlights the need for each system parameter to be within a certain range of values in order for oscillations to be observed. For the cGENIE model, this implies that oscillations could easily be missed if, for example, the wind speed over the southern channel is too large (see \ref{app:wind} for a discussion on the role of wind speed).

The results presented here provide a possible explanation for observations in cGENIE experiments using realistic continental configurations of the geological past, where stable oscillations coincide with configurations having a large open ocean pole \cite{pohl22}. These cGENIE oscillations, which occur on the millennial timescale, could shed light on periods of high variability observed in high-resolution proxy data \cite{Dahl2019}.

Despite some evidence of stable oscillatory behaviour in GCMs \cite{valdes21}, the extended run-time of such models makes them ill-suited for the explorations through parameter space required when investigating possible causes and attributing factors of oscillatory behaviour.
cGENIE, on the other hand, represents a very practical compromise for establishing connections between interesting dynamics and fundamental mechanisms that drive them. As such, this study is conducted in a similar vein to some previous studies that have conducted parameter sweeps in order to identify possible bifurcations and regions of qualitatively different climate states (for example, see Refs.~\cite{crichton2021,marsh2013,marsh2004}).

The way that the stable oscillations presented here are generated by delayed feedback (in turn, by continental configuration) is different to previous studies of ocean models. In some past studies convection and associated deep-water flushes have been the driving mechanism for oscillations \cite{Welander1982,ColindeVerdiere2007}.
In other studies sea-ice dynamics are required as part of a feedback mechanism that generates oscillations \cite{li2019,roberts2017,wang2006}.
One logical next step is to investigate how the delayed feedback interacts with these other sources of oscillations. There is evidence from cGENIE simulations (not shown) of quasi-periodic behaviour, which could be the result of an interaction between multiple mechanisms for oscillations.

The simple box model with delayed feedback does not reproduce all dynamical features observed in the cGENIE experiments, such as the oscillations shown as blue filled circles in Fig.~\ref{fig:amp_large}(a). Furthermore, although the simple model captures the millennial timescale of the stable oscillations, the period of the oscillations observed in cGENIE can easily exceed 10 thousand years (for example, the large-amplitude oscillations around 9 PAL in Fig.~\ref{fig:genie}(c)), which is not found in the simple model. A further modelling possibility that would effect the timescales in the box model includes adding more boxes to the box model. 
Finally, deep-water flushes could be captured by a salinity delay due to a deep-ocean adjustment to convective events, which would further increase the delay time of the model.

\begin{figure}[t]
  \centering
  \vspace*{1mm}
  \includegraphics[width=\columnwidth]{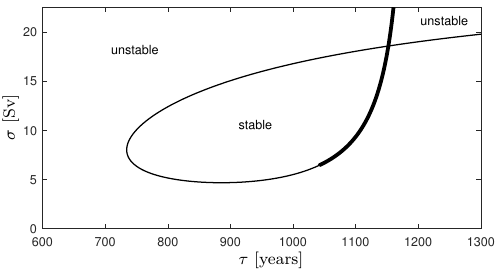}
  \caption{Stability diagram of the steady-state solution at $F_1=-0.22 \text{~Sv}$ in the $(\tau,\sigma)$-plane. Non-bold/bold curves show subcritical/supercritical Hopf bifurcations. 
}
  \label{fig:hopf_sigma_tau}
\end{figure}

Results of our box model may also have implications for our understanding of future changes in the global ocean circulation. It has been reported that part of the Antarctic Circumpolar Current (ACC) is accelerating due to anthropogenic ocean warming \cite{shi21}. This would correspond to a shortened delay time in the feedback of model~(\ref{eq:titz_delay}). Figure~\ref{fig:hopf_sigma_tau} shows a curve of Hopf bifurcations that demonstrates the sensitivity of an example ${m>0}$ solution ($F_1=-0.22 \text{~Sv}$) to changes in the feedback. Interestingly, it is stable only for intermediate values of feedback strength and delay time. In both scenarios of increasing and decreasing delay time, the steady-state solution loses stability at a Hopf bifurcation. For the parameters used in Fig.~\ref{fig:hopf_sigma_tau}, there are segments of supercritical Hopf bifurcations (bold), which give rise to stable oscillations, as well as, subcritical Hopf bifurcations (non-bold). The latter correspond to so-called ``dangerous'' bifurcations \cite{thompson11}, where the system jumps suddenly and irreversibly to a remote state with an alternative circulation pattern. Therefore, the dynamics discussed in this paper could also be relevant in a modern context. Of course, the Drake Passage of modern continental configuration is smaller than the passage of the Drake World setup used here (Fig.~\ref{fig:DW}(a)). Nonetheless, it is worthwhile to consider the possibility that changes in zonal flow, such as that of the ACC, could play a part in destabilising the meridional overturning circulation of the Atlantic Ocean.

\section*{Data availability}

The code for the version of the ‘muffin’ release of the cGENIE Earth system model used in this paper, is tagged as v0.9.58, and is assigned a DOI: 10.5281/zenodo.14606515.
Configuration files for the specific experiments presented in the paper can be found in the
directory: genie-userconfigs/PUBS/submitted/
Keane\_et\_al.PhysicaD.2025. Details of the experiments, plus the command line needed to run each one, are given in the readme.txt file in that directory. All other configuration files and boundary conditions are provided as part of
the code release.
A manual detailing code installation, basic model configuration, tutorials covering various aspects of model configuration, experimental design, and output, plus the processing of
results, is assigned a DOI: 10.5281/zenodo.1407657. 

We provide code for the numerical continuation of model~(\ref{eq:titz_delay}): https://github.com/andrewkeane/kpdr\_oscillations
\noindent It demonstrates the calculation of curves of steady-state solutions and Hopf bifurcations using the Matlab-based software DDE-Biftool.

\section*{Acknowledgements}

This is a contribution to UMR 6282 Biogeosciences Team SEDS and UNESCO project IGCP 735 ``Rocks and the Rise of Ordovician Life'' (Rocks n' ROL). AP acknowledges the support of the French Agence Nationale de la Recherche (ANR) under references ANR-22-CE01-0003 (project ECO-BOOST) and ANR-23-CE01-0003 (project CYCLO-SED) and the Leverhulme Trust for supporting the project `Earth System dynamics at the dawn of the animal-rich biosphere' (RPG-2022-233). Calculations were performed using HPC resources from DNUM CCUB (Centre de Calcul de l'Universit\'e de Bourgogne).
The project has received funding from the ERC, project 101055096 (HD). 
AR acknowledges support from NSF grant  EAR-2121165, as well as from the Heising-Simons Foundation.
This work benefited from a research visit of AK to University of California-Riverside funded by a Charlemont grant from the Royal Irish Academy.

We are grateful to two anonymous reviewers for their feedback that helped to improve the quality and presentation of the paper.

\appendix
\section{Model equations for negative $m$}
\label{app:neg_m}
In the original model, without the delayed feedback term and as described in Ref.~\cite{titz02}, one set of equations is applied when the ocean flow is clockwise, while another set of equations is applied when the flow is counter-clockwise. The direction of the flow is determined by the model variable $m$, as defined in Eq.~\ref{eq:m} above. Since the oscillatory dynamics that we are interested in occur for positive $m$, we only provide the model equations for the positive $m$ case in the main text above. When $m$ is negative, the meridional flux between boxes (shown in Fig.~\ref{fig:boxes}) is reversed and Eq.~\ref{eq:titz_delay} becomes:
\begin{align}
\label{eq:titz_delay_negative_m}
V\dot{S_1} &= S_0 F_1 - m \left(S_3 - S_1\right) + \sigma \left(S_1^{\tau} - S_1\right)\\ \nonumber
V\dot{S_2} & =  -S_0 F_2 - m \left(S_1 - S_2\right).
\end{align}

\section{Salt conservation approximation}
\label{app:salt}

The salt conservation rule used in the above analysis 
\begin{equation}
\label{eq:salt}
    S_3(t) = 3S_0 - S_1(t) - S_2(t)
\end{equation}
takes the same form as in previous studies of ocean box models \cite{titz02, titz02b} and has the significant advantage of reducing the model to a relatively simple set of two differential equations (i.e. Eqs.~\eqref{eq:titz_delay}). 
This rule can be interpreted as saying that any salt from the total $3VS_0$ that is not in box~1 or box~2 must be in box~3.
It can be derived by summing the equations of the full three box model, without the delay term,
\begin{align*}
    V\dot{S_1} &= S_0 F_1 + m \left(S_2 - S_1\right) \\ \nonumber
    V\dot{S_2} & =  -S_0 F_2 + m \left(S_3 - S_2\right)\\ \nonumber
    V\dot{S_3} & = S_0 (F_2 - F_1) + m \left(S_1 - S_3\right)
\end{align*}
to find
\begin{equation*}
    \dot{S_1} + \dot{S_2} + \dot{S_3} = 0
\end{equation*}
then integrating from $t_0$ to $t$
\begin{align*}
    & S_1(t) + S_2(t) + S_3(t)\nonumber \\ 
    & - S_1(t_0) - S_2(t_0) - S_3(t_0) = 0\\
    \Rightarrow & S_1(t) + S_2(t) + S_3(t) - 3S_0 = 0\\
    \Rightarrow & S_3(t) = 3S_0 - S_1(t) - S_2(t).
\end{align*}

However, once the delay term is included in the model, this calculation only holds when $S^\tau_1=S_1$. Otherwise, Eq.~\eqref{eq:salt} serves only as an approximation.
To derive a more precise salt conservation rule, we begin with the three box model with the delay term
\begin{align}
    \label{eq:full3_delay}
    V\dot{S_1} &= S_0 F_1 + m \left(S_2 - S_1\right) + \sigma \left(S_1^{\tau} - S_1\right)\\ \nonumber
    V\dot{S_2} & =  -S_0 F_2 + m \left(S_3 - S_2\right)\\ \nonumber
    V\dot{S_3} & = S_0 (F_2 - F_1) + m \left(S_1 - S_3\right).
\end{align}

One might be inclined to compare our analysis in section~\ref{section:oscillations} with an analysis of model~\eqref{eq:full3_delay}. However, numerical continuation of model~\eqref{eq:full3_delay} proves awkward, since for $S^\tau_1=S_1$ (i.e. for all fixed points) the Jacobian of the system becomes singular and numerical continuation cannot be applied. Instead, we derive an updated salt conservation rule by summing Eqs.~\eqref{eq:full3_delay} to get
\begin{align}
\label{eq:salt_new1}
    \dot{S_1} + \dot{S_2} + \dot{S_3} &= -\frac{\sigma}{V} \left(S_1 - S_1^{\tau}\right).
\end{align}
In order to find a salt conservation rule that is relatively easy to work with (e.g. without any derivatives or integrals) we apply the fundamental theorem of calculus to get
\begin{align}
\label{eq:salt_new2}
    S_1(t) - S_1(t-\tau) = \frac{d}{dt} \int^t_{t-\tau} S_1(t') dt' 
\end{align}
and approximate the integral using the composite trapezoidal rule

\begin{align}
\label{eq:salt_new3}
    \int^t_{t-\tau} S_1(t') dt' & = \frac{0.1\tau}{2}( S_1(t) + 2S_1(t-0.1\tau)\nonumber \\ 
    & + 2S_1(t-0.2\tau) + 2S_1(t-0.3\tau)\nonumber \\ 
    & + 2S_1(t-0.4\tau) + 2S_1(t-0.5\tau)\nonumber \\ 
    & + 2S_1(t-0.6\tau) + 2S_1(t-0.7\tau)\nonumber \\ 
    & + 2S_1(t-0.8\tau) + 2S_1(t-0.9\tau)\nonumber \\ 
    & + S_1(t-\tau)).
\end{align}
Combining Eqs.~\eqref{eq:salt_new1}--\eqref{eq:salt_new3} gives
\begin{align*}
    & \dot{S_1}(t) + \dot{S_2}(t) + \dot{S_3}(t) = -\frac{\sigma}{V} \left(S_1(t) - S_1(t-\tau)\right)\nonumber \\
    &= -\frac{\sigma}{V}\frac{0.1\tau}{2}(\dot S_1(t) + 2\dot S_1(t-0.1\tau) + 2\dot S_1(t-0.2\tau)\nonumber \\ 
    & + 2\dot S_1(t-0.3\tau) + 2\dot S_1(t-0.4\tau) + 2\dot S_1(t-0.5\tau)\nonumber \\ 
    & + 2\dot S_1(t-0.6\tau) + 2\dot S_1(t-0.7\tau) + 2\dot S_1(t-0.8\tau)\nonumber \\ 
    & + 2\dot S_1(t-0.9\tau) + \dot S_1(t-\tau)).
\end{align*}
This leads to the updated salt conservation rule
\begin{align}
\label{eq:salt_new}
    S_3(t) &= 3S_0 - S_1(t) - S_2(t) -\frac{\sigma}{V}\frac{0.1\tau}{2}(S_1(t)\nonumber \\ 
    & + 2S_1(t-0.1\tau) + 2S_1(t-0.2\tau)\nonumber \\ 
    & + 2S_1(t-0.3\tau) + 2S_1(t-0.4\tau)\nonumber \\ 
    & + 2S_1(t-0.5\tau) + 2S_1(t-0.6\tau)\nonumber \\ 
    & + 2S_1(t-0.7\tau) + 2S_1(t-0.8\tau)\nonumber \\ 
    & + 2S_1(t-0.9\tau) + S_1(t-\tau)).
\end{align}

Equation~\eqref{eq:salt_new} may appear cumbersome, but replacing Eq.~\eqref{eq:salt} with Eq.~\eqref{eq:salt_new} in DDE-Biftool is very simple. It only takes several lines of code to define the additional delay times required for the new salt conservation rule.

\begin{figure}[t]
  \centering
  \vspace*{1mm}
  \includegraphics[width=\columnwidth]{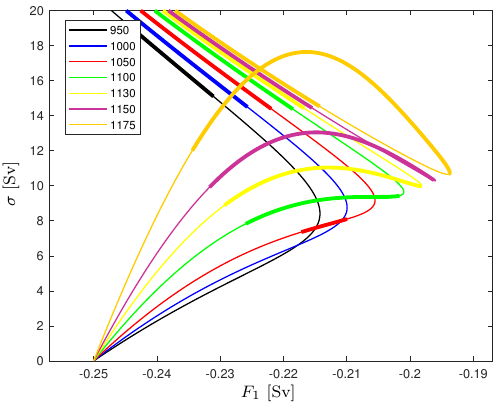}
  \caption{Same as Fig.~\ref{fig:hopf} with more realistic salt conservation rule (Eq.~\eqref{eq:salt_new}).
  }
  \label{fig:hopf_newsalt}
\end{figure}

We analyse the effect that this new rule has on the results presented above. To keep the comparison brief we show the equivalent of Fig.~\ref{fig:hopf} with the new salt rule, since this figure shows the main effect of the delayed feedback stabilising periodic orbits in a concise way.

Figure~\ref{fig:hopf_newsalt} displays curves of Hopf bifurcations in the $(F_1,\sigma)$-plane for different values of $\tau$ using model~\eqref{eq:titz_delay} with Eq.~\eqref{eq:salt_new}. Again, thin and thick curves indicate subcritical and supercritical bifurcations. We observe that, although there are some quantitative differences between Figs.~\ref{fig:hopf} and~\ref{fig:hopf_newsalt}, the updated salt conservation rule does not impede the role of the delayed feedback in stabilising periodic orbits or changing the criticality of Hopf bifurcations from subcritical to supercritical.

\section{Steady states of model~(\ref{eq:titz_delay})}
\label{app:steady}

Here we include calculations on the linear stability analysis of steady-state solutions to model~(\ref{eq:titz_delay}). Due to the nature of DDEs, the characteristic equation is transcendental, which, in this case, limits the analysis to finding a pair of expressions describing the Hopf bifurcations.

First, we nondimensionalise the model by rescaling the dependent variables
\begin{align*}
    \hat S_1 &= \frac{S_1}{S_0},\\
    \hat S_2 &= \frac{S_2}{S_0},\\
    \hat S_3 &= 3 - \hat S_1 - \hat S_2.
\end{align*}

We rescale time, $t$, such that
\begin{align*}
    \hat t &= \frac{kt}{V},\\
    \hat S_1^{\hat\tau} &= \hat S_1(\hat t-\hat\tau) \text{ with } \hat\tau = \frac{k\tau}{V}.
\end{align*}

The parameters are rescaled as:
\begin{align*}
    \hat F_1 &= \frac{F_1}{k},\\
    \hat F_2 &= \frac{F_2}{k},\\
    \hat \sigma &= \frac{\sigma}{k}, \\
    \hat \beta &= \beta S_0,\\
    \hat T^* &= \alpha T^*.
\end{align*}

The resulting nondimensionalised model is
\begin{align*}
\frac{d\hat S_1}{d\hat t} &= \hat F_1 + \left[\hat\beta\left(\hat S_2 - \hat S_1\right) - \hat T^*\right]\left(\hat S_2 - \hat S_1\right) \\
&+ \hat\sigma\left(\hat S_1^{\hat\tau} - \hat S_1\right)\\
\frac{d\hat S_2}{d\hat t} & =  -\hat F_2 + \left[\hat\beta\left(\hat S_2 - \hat S_1\right) - \hat T^*\right]\left(\hat S_3 - \hat S_2\right).
\end{align*}

To make further calculations easier, we rewrite the variables as
\begin{align*}
    \hat S_1 &= S_1, \\
    \hat S_{21} &= \hat S_2 - \hat S_1, \\
    \hat S_3 &= 3 - 2\hat S_1 - \hat S_{21}.
\end{align*}

By letting $\hat F_{21} = \hat F_2 + \hat F_1$, the model can be written as
\begin{align*}
\frac{d\hat S_1}{d\hat t} &= \hat F_1 + \left(\hat\beta\hat S_{21} - \hat T^*\right)\hat S_{21} + \hat\sigma\left(\hat S_1^{\hat\tau} - \hat S_1\right)\\ 
\frac{d\hat S_{21}}{d\hat t} &= -\hat F_{21} + 3\left(\hat\beta\hat S_{21} - \hat T^*\right)\left(1 - \hat S_1 - \hat S_{21}\right) \\
&- \hat\sigma\left(\hat S_1^{\hat\tau} - \hat S_1\right).
\end{align*}

For simplicity, we let $\hat T^*=0$, as is done for the results in the main text, and solve for steady-state solutions, $(\hat S_1^*,\hat S_{21}^*)$, with $\hat S_1^{\hat\tau} = \hat S_1$:
\begin{align*}
\hat S_{21}^* &= \pm \sqrt{\frac{-\hat F_1}{\hat\beta}},\\ 
\hat S_1^* &= 1 - \hat S_{21}^* - \frac{\hat F_{21}}{3\hat\beta\hat S_{21}^*}.
\end{align*}

To consider the stability of the steady states, we linearise the system around $(\hat S_1^*,\hat S_{21}^*)$:

\begin{align*}
\begin{pmatrix}
\frac{d\hat{S}_1}{d\hat t} \\
\frac{d\hat{S}_{21}}{d\hat t}
\end{pmatrix} &=
\begin{pmatrix}
-\hat\sigma & 2\hat\beta\hat{S}_{21}^*\\
-3\hat\beta\hat S_{21}^* +\hat\sigma & 3\hat\beta - 3\hat\beta\hat S_1^*- 6\hat\beta\hat{S}_{21}^*
\end{pmatrix}
\begin{pmatrix}
\hat{S}_1 \\
\hat{S}_{21}
\end{pmatrix} \\
&+
\begin{pmatrix}
\hat\sigma & 0\\
-\hat\sigma & 0
\end{pmatrix}
\begin{pmatrix}
\hat{S}_1^{\hat\tau} \\
\hat{S}_{21}^{\hat\tau}
\end{pmatrix}.
\end{align*}

As derived in Ref.~\cite{smith2011introduction}, the characteristic equation of a linear system of $n$ DDEs
\begin{equation*}
    \frac{d}{dt}\mathbf{x}(t) = \mathbf{A}\mathbf{x}(t) + \mathbf{B}\mathbf{x}(t-\tau), \quad \mathbf{x}\in\mathbb{R}^n,
\end{equation*}
is given by
\begin{equation*}
    \det{\{\lambda\mathbb{I} - \mathbf{A} - e^{-\lambda\tau} \mathbf{B}\}}.
\end{equation*}

By evaluating at the steady state, noting that the model is designed for $\hat{S}_{21}>0$, and letting
\begin{align*}
    a &= \sqrt{-\hat\beta\hat F_1}, \\
    b &= \hat\beta\hat F_{2},
\end{align*}

the characteristic equation becomes
\begin{equation*}
\det{\begin{pmatrix}
\lambda + \hat\sigma - \hat\sigma e^{-\lambda\hat\tau}  & -2a\\
3a - \hat\sigma + \hat\sigma e^{-\lambda\hat\tau} & \lambda + 4a - \frac{b}{a} 
\end{pmatrix}}=0.
\end{equation*}

We are specifically interested in Hopf bifurcations. Therefore, we set $\lambda = i\hat\omega$, where $\hat\omega$ corresponds to the nondimensionalised frequency of oscillations around the steady state. By separating the resulting expression into real and imaginary parts, we find the following pair of expressions:
\begin{align*}
0 &= 6a^3 + 2a^2\hat\sigma(1 - \cos(\hat\omega\hat\tau)) -a\hat\omega(\hat\omega + \hat\sigma\sin(\hat\omega\hat\tau)) \\
&- b\hat\sigma(1 - \cos(\hat\omega\hat\tau)), \\
0 &= 2a^2(2\hat\omega + \hat\sigma\sin(\hat\omega\hat\tau)) + a\hat\omega\hat\sigma(1 - \cos(\hat\omega\hat\tau)) \\
&- b(\hat\omega + \hat\sigma\sin(\hat\omega\hat\tau)).
\end{align*}

These expressions can be verified with the Hopf bifurcations found numerically in the main text. However, using them for further analysis is impractical and at this point it makes sense to rely on numerical techniques.

\section{Oscillation dependence on wind speed}
\label{app:wind}

The wind speed in cGENIE is physically related to the delay time in model~(\ref{eq:titz_delay}). In particular, the zonal flow across the southern channel (i.e. the circumpolar current) is primarily wind-driven. Therefore, the wind speed across the latitudes of the southern channel in cGENIE strongly effects the average zonal flow velocity of the circumpolar current. On the other hand, the delay time in model~(\ref{eq:titz_delay}) represents the transportation time of water around the southern pole, averaged across all depths and latitudes of the southern channel. For each latitude, the transportation time is inversely proportional to the average zonal flow velocity. Hence, there is an inverse relationship between the wind speeds across the southern channel in cGENIE and the average transportation time, represented by the delay time in model~(\ref{eq:titz_delay}).

\begin{figure}[t]
  \centering
  \vspace*{1mm}
  \includegraphics[width=\columnwidth]{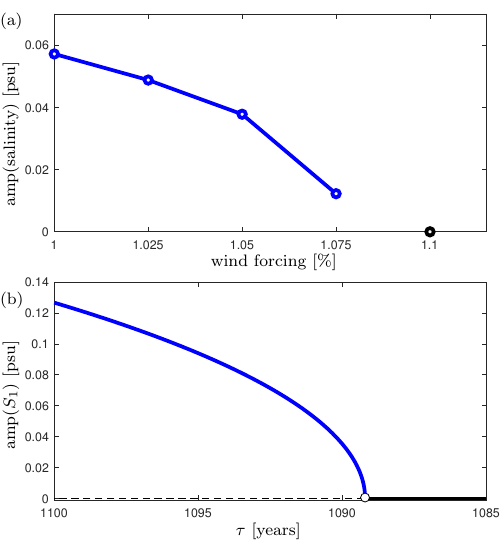}
  \caption{Amplitudes of oscillations with steady-state and periodic solutions shown in black and blue, respectively. (a) cGENIE experiments for 6 PAL with varying wind force strength across the southern channel as a factor of the strength used throughout the paper. (b) Model~(\ref{eq:titz_delay}) solutions while varying $\tau$ with $\sigma=9.5 \text{~Sv}$ and $F_1=0.21\text{~Sv}$. White filled circle represents Hopf bifurcations.
  }
  \label{fig:wind}
\end{figure}

To highlight the sensitivity of the stable oscillations to model parameters of cGENIE, we show that the oscillations begin to disappear as the wind speed over the southern channel is slightly increased. 
Figure~\ref{fig:wind}(a) shows the amplitude of the stable oscillations in the cGENIE model, in terms of mean salinity (across the southern channel), of the solution at 6 PAL. The horizontal axis represents the wind speed, across the southern channel, as a factor of the wind speed used in the rest of the paper, as illustrated in Fig.~\ref{fig:DW}(b). In other words, the solution with wind speed equal to one corresponds to the same 6 PAL solution shown in Figs.~\ref{fig:genie} and~\ref{fig:amp_small}. As the wind speed is increased we observe that the amplitude of the oscillations decreases and appears to be approaching zero.

Since the circulation around the pole is wind-driven, the wind speed determines the velocity of the flow, and therefore determines the delay time associated with the delayed feedback in the box model. To demonstrate this relationship between wind speed of the cGENIE model and the delay time of the box model, we consider the box model solution shown in Fig.~\ref{fig:amp_small}(b) at $F_1=-0.21$ that approximately corresponds to the cGENIE oscillation at 6 PAL. In Fig.~\ref{fig:wind}(b) we observe that as the delay time is decreased, the amplitude of the oscillations decreases until reaching zero at a Hopf bifurcation.



\bibliographystyle{ieeetr} 
\bibliography{kpdr_oscillations_refs}

\end{document}